\documentclass[11pt,reqno]{amsart}

\usepackage{amsmath}
\usepackage{amsfonts}
\usepackage{amscd}
\usepackage{enumerate}
\usepackage{fancyheadings}
\usepackage{bbold}

\usepackage{amsmath}
\usepackage{amsthm}
\usepackage{amssymb}
\usepackage{amsfonts}
\usepackage{epsfig}
\usepackage{epsf}
\newtheorem{theorem}{Theorem}[section]

\newtheorem{example}{Example}[section]

\newtheorem*{remark}{Remark}

\newtheorem{lemma}{Lemma}[section]

\usepackage{latexsym}
\usepackage{epsfig}

\newcommand{\baseRing}[1]{\ensuremath{\mathbb{#1}}}
\newcommand{\C}{\baseRing{C}}

\newcommand{\g}{\mathfrak{g}}
\newcommand{\h}{\mathfrak{h}}
\renewcommand{\b}{\mathfrak{b}}

\renewcommand{\phi}{\varphi}
\allowdisplaybreaks
\begin{document}

\title{Super Solutions of the Dynamical Yang-Baxter Equation}
\author{Gizem Karaali}
\address{Department of Mathematics, University of California, Santa Barbara, Ca
93106}
\email{gizem@math.ucsb.edu}

\begin{abstract}
A super dynamical r-matrix $r$ satisfies the \emph{zero weight
condition} if:
$$ [h\otimes 1 + 1 \otimes h, r(\lambda)] = 0 \textmd{ for all }
h \in \h, \lambda \in \h^*.$$
\noindent
In this paper we classify super dynamical $r-$matrices with zero weight, thus extending the results
of \cite{EV} to the graded case. 
\end{abstract}

\maketitle

\section{Introduction}
\label{SectionIntroduction}

A complete classification of the nonskewsymmetric solutions of the classical
Yang-Baxter equation
exists in the case when the underlying structure is a simple Lie algebra; see
\cite{BD1} and
\cite{BD2} for the original proofs by Belavin and Drinfeld, and \cite{ES} for a
more pedagogical
exposition. A similar construction, with natural modifications, works in the
super case as well;
see \cite{Kar1}. However, it turns out that this may not be easily modified
into a full
classification result; see \cite{Kar2} for an explicit construction and
detailed study of a
counterexample. 

It is well-known that solutions of the classical Yang-Baxter equation on a Lie
algebra give us the
semiclassical limits of quantizations on the associated Lie group. In
\cite{ESS}, Etingof, Schedler
and Schiffmann have explicitly constructed quantizations associated to all
solutions coming from
the Belavin-Drinfeld result. Their method in fact works for all dynamical $r-$matrices, 
i.e. the solutions of the more general dynamical Yang-Baxter equation.

The purpose of this paper is to begin a study of the super analog of the theory of
dynamical $r-$matrices.
Ultimately our goal is a full theory of quantum groups in the super setting. We
expect that
understanding the super solutions of the dynamical Yang-Baxter equation will
provide us with
valuable insight, and hence help us extend or modify the quantization result
cited above 
to obtain a graded analogue.\footnotemark \, A possible classification result
in this more
general setting
of dynamical $r-$matrices may also clarify the so far exceptional case of the
classical $r-$matrices
mentioned before. 

\footnotetext{As pointed out by P. Etingof, a different method of quantization, by Etingof and
Kazhdan, has
been generalized to
the super setting, see \cite{Ge}. However, this is a less constructive
technique, and does not
fully answer our questions.}


\section{Dynamical \textit{r}-matrices in the Super Setting}
\label{SectionSuperDynamicalrmatrices}

\subsection{Definitions}
\label{DefinitionsSuperDynamical}

Let $\g$ be a simple Lie superalgebra with non-degenerate Killing
form ${(\cdot \;,\cdot)}$. Let ${\h \subset \g}$ be a Cartan
subsuperalgebra, and let $\Delta \subset \h^*$ be the set of
roots associated to $\h$. Fix a set of simple roots $\Gamma$ or
equivalently a Borel $\b$. The \emph{classical dynamical
Yang-Baxter equation} for a meromorphic function ${r:\h^*
\rightarrow \g \otimes \g}$ will be: 
\begin{equation}
\label{SDYBequation}
Alt_s(dr) + 
[r^{12},r^{13}] + [r^{12},r^{23}] + [r^{13},r^{23}] = 0
\end{equation}
\noindent
The differential of $r$ will be defined as
above as:
$$ \begin{matrix}
dr &:& \h^* &\longrightarrow& \g \otimes \g \otimes \g  \\
&& \lambda &\longmapsto& \sum_i x_i \otimes \frac{\partial
r}{\partial x_i} (\lambda)  
\end{matrix}. $$
\noindent
Here $\{x_i\}$ is a basis for $\h$ so all $x_i$ are even.
Recall that ${Alt_s : \g^{\otimes 3} \rightarrow
\g^{\otimes 3}}$ is given on homogeneous elements by: 
$$ Alt_s(a\otimes b \otimes c) = a \otimes b \otimes c  +
(-1)^{|a|(|b|+|c|)} b \otimes c \otimes a + 
(-1)^{|c|(|a|+|b|)} c \otimes a \otimes b, $$
\noindent
In view of all this, we can see that for $r = \sum_i
{R_i}_{(1)} \otimes {R_i}_{(2)}$:
\begin{align*}
Alt_s(dr) =& 
\sum_i x_i^{(1)} \left(\frac{\partial r }{\partial x_i}
\right)^{(23)} +
\sum_i x_i^{(2)} \left(\frac{\partial r }{\partial x_i}
\right)^{(31)}\\ \\
+& 
\sum_i (-1)^{|{R_i}_{(1)}| |{R_i}_{(2)}|} x_i^{(3)} \left(
\frac{\partial r}{\partial x_i}\right)^{(12)}.
\end{align*}

We will say that a meromorphic function ${r : \h^* \rightarrow \g
\otimes \g}$ is a \emph{super dynamical r-matrix with coupling
constant} $\epsilon$ if it is a solution to Equation
\ref{SDYBequation} and satisfies the generalized unitarity
condition:
\begin{equation}
\label{dynamicalunitarity} 
r(\lambda) + T_s(r)(\lambda) = \epsilon \Omega,
\end{equation}
\noindent
where $\Omega$ is the Casimir element, i.e. the element of ${\g
\otimes \g}$ corresponding to the Killing form. Here, ${T_s : V \otimes V  \rightarrow V
\otimes V}$ is the
\emph{super twist map} defined on the homogeneous elements of a given super
vector space $V$ as  
$${T_s (a\otimes b) = (-1)^{|a||b|}b\otimes a}.$$

\begin{remark}
The above definitions are the natural super analogues of the non-graded terminology used in the
study of the dynamical Yang-Baxter equation,
see \cite{Et} for a survey on the non-graded theory. 
\end{remark}

\subsection{Super Dynamical r-matrices with Zero Weight}
\label{SuperZeroWeight}

A super dynamical r-matrix $r$ satisfies the \emph{zero weight
condition} if:
$$ [h\otimes 1 + 1 \otimes h, r(\lambda)] = 0 \textmd{ for all }
h \in \h, \lambda \in \h^*.$$
\noindent
In the next two sections we will prove the following two
statements:

\vspace{0.1in}
\begin{theorem}
\label{0couple0weighttheorem}
(1) Let $X$ be a subset of the set of roots $\Delta$ of a
simple Lie superalgebra $\g$ with non-degenerate Killing form ${(
\cdot \; , \cdot )}$ such that:

(a) If $\alpha, \beta \in X$ and $\alpha + \beta$ is a root,
then $\alpha + \beta \in X$, and

(b) If $\alpha \in X$, then $-\alpha \in X.$

\noindent
Let $\nu \in \h^*$, and let $D = \sum_{i<j} D_{ij} dx_i \wedge
dx_j$ be a closed meromorphic $2-$form on $\h^*$. If we set
$D_{ij} = -D_{ji}$ for $i \ge j$, then the meromorphic
function:   
$$ r(\lambda) =
\sum_{i,j=1}^N D_{ij}(\lambda) x_i \otimes x_j  +  \sum_{\alpha
\in X} \frac{(-1)^{|\alpha|}(e_{\alpha},e_{-\alpha})}{(\alpha,
\lambda-\nu)} e_{\alpha} \otimes e_{-\alpha} $$ 
\noindent
is a super dynamical r-matrix with zero weight and zero
coupling constant.

\noindent
(2) Any super dynamical r-matrix with zero weight and zero
coupling constant is of this form. 
\end{theorem}

\vspace{0.1in}
\begin{theorem}
\label{0weighttheorem}
(1) Let $\g$ be a simple Lie superalgebra with non-degenerate
Killing form ${(\cdot \; , \cdot )}$. Let $X$ be a subset of the set of
roots $\Delta$ of $\g$ satisfying conditions $(a)$ and $(b)$ of
Theorem \ref{0couple0weighttheorem}.  
Pick $\nu \in \h^*$, and define:
$$ \phi_{\alpha} = \left\{ \begin{matrix}
\left( \epsilon / 2\right) \emph{coth}\left(
(-1)^{|\alpha|}(e_{\alpha},e_{-\alpha})
\left(\epsilon / 2\right) (\alpha,\lambda - \nu) \right) 
& \textmd{   if  } \alpha \in X \\
\left(\pm\epsilon / 2\right) 
& \textmd{   if  } \alpha \not \in X, 
\textmd{  negative} \\ 
\mp(-1)^{|\alpha|}\left(\epsilon / 2 \right) 
& \textmd{   if  } \alpha \not \in X,
\textmd{  positive}  
\end{matrix} \right. $$
\noindent
Let $D = \sum_{i<j} D_{ij} dx_i \wedge
dx_j$ be a closed meromorphic $2-$form on $\h^*$. If we set
$D_{ij} = -D_{ji}$ for $i \ge j$, then the meromorphic
function:    
$$ r(\lambda) =
\sum_{i,j=1}^N D_{ij}(\lambda) x_i \otimes x_j +
\frac{\epsilon}{2}\Omega + \sum_{\alpha \in \Delta}
\phi_{\alpha}
e_{\alpha} \otimes e_{-\alpha} $$  
\noindent
is a super dynamical r-matrix with zero weight and nonzero
coupling constant $\epsilon$. 

\noindent
(2) Any super dynamical r-matrix with zero weight and nonzero
coupling constant $\epsilon$ is of this form. 
\end{theorem}

\vspace{0.1in}
\begin{remark}
Note that if we take the limit as
$\epsilon \rightarrow 0,$ the above expression reduces to the
expression of Theorem \ref{0couple0weighttheorem}.
\end{remark}

\subsection{Proof of Theorem \ref{0couple0weighttheorem}}
\label{Proofof00Thm}

For any positive root $\alpha$ fix $e_{\alpha} \in \g_{\alpha}$
and pick ${e_{-\alpha} \in \g_{-\alpha}}$ dual to $e_{\alpha}$
i.e. 
$$ (e_{\alpha}, e_{-\alpha}) = 1 \textmd{ for all } \alpha \in
\Delta^+.$$

We introduce the following notation:
$$A_{\alpha} = \left\{ \begin{array}{cl} 
(-1)^{|\alpha|} & \textmd{if } \alpha  \textmd{ is positive} \\
1 & \textmd{if } \alpha \textmd{ is negative}
\end{array} \right. $$
\noindent
Note that $A_{-\alpha} = (-1)^{|\alpha|}A_{\alpha}$.
We can use $A_{\alpha}$ for instance to write the duals of our
basis vectors in terms of one another:
$$ e_{\alpha}^* = A_{-\alpha} e_{-\alpha}$$
\noindent
or equivalently:
$$ (e_{\alpha},e_{-\alpha}) = A_{-\alpha}. $$

Let $r : \h^* \rightarrow \g \otimes \g$ be a super dynamical
r-matrix with zero coupling constant. Then the zero weight
condition on $r$ implies that $r$ has to be of the form:
$$ r(\lambda) = \sum_{i,j} D_{ij}(\lambda) x_i \otimes x_j +
\sum_{\alpha \in \Delta} \phi_{\alpha} e_{\alpha} \otimes
e_{-\alpha}, $$
\noindent 
where $D_{ij}, \phi_{\alpha}$ are suitable scalar meromorphic
functions such that:
$$ D_{ij}(\lambda) = - D_{ji}(\lambda)  
\; \; \; \; \;
\textmd{and} 
\; \; \; \; \;
\phi_{-\alpha} = -(-1)^{|\alpha|} \phi_{\alpha}$$
\noindent
One can easily check that there can be no terms mixing the
Cartan part with the non-Cartan part, and the conditions on the
$D_{ij}$ and the $\phi_{\alpha}$ follow from the zero coupling
constant.

Equation \ref{SDYBequation} is an equation in ${\g \otimes \g
\otimes \g}$. The zero coupling constant implies that the left
hand side of the equation is skew-symmetric with respect to
signed permutations of factors: 
\begin{eqnarray*}
(12)_s (a \otimes b \otimes c) &=& (-1)^{|a||b|} b \otimes a
\otimes c \\
(13)_s (a \otimes b \otimes c) &=& (-1)^{|a||b| + |a||c| + |b||c|} c
\otimes b \otimes a
\\
(23)_s (a \otimes b \otimes c) &=& (-1)^{|b||c|} a \otimes c
\otimes b 
\end{eqnarray*}

Therefore in order to solve Equation \ref{SDYBequation}, it is
enough to look at its ${\h \otimes \h \otimes \h}$, ${\h \otimes
\g_{\alpha} \otimes \g_{-\alpha},}$ and ${\g_{\alpha} \otimes
\g_{\beta} \otimes \g_{-\alpha -\beta}}$ parts.

The $\h \otimes \h \otimes \h$ part is:
\begin{eqnarray*}
Alt_s(dr) \left|_{\h \otimes \h \otimes \h}\right. 
&=& 
\sum_{i,j,k} 
x_i^{(1)} \frac{\partial D_{jk}x_jx_k^{(23)}}{\partial x_i}  +  
\sum_{i,j,k} 
x_i^{(2)} \frac{\partial D_{jk}x_jx_k^{(31)}}{\partial x_i} \\
&&+  
\sum_{i,j,k} 
x_i^{(3)} \frac{\partial D_{jk}x_jx_k^{(12)}}{\partial x_i}
\\ \\
&=& \sum_{i,j,k}  
\frac{\partial D_{jk}}{\partial x_i} x_i \otimes x_j \otimes x_k 
+ \sum_{i,j,k} 
\frac{\partial D_{jk}}{\partial x_i} x_k \otimes x_i \otimes x_j \\
&&+ \sum_{i,j,k} 
\frac{\partial D_{jk}}{\partial x_i} x_j \otimes x_k \otimes x_i
\\ \\
&=& \sum_{i,j,k} 
\left( \frac{\partial D_{ij}}{x_k} + \frac{\partial D_{jk}}{x_i} 
+ \frac{\partial D_{ki}}{x_j} \right) 
x_i \otimes x_j \otimes x_k,  
\end{eqnarray*}
\noindent
and it vanishes if and only if $D = {\sum_{i < j} D_{ij} dx_i
\wedge dx_j}$ is a closed $2-$form.

To find the $\h \otimes \g_{\alpha} \otimes \g_{-\alpha}$ part,
we first look at $[[r,r]]$:
\begin{eqnarray*}
[r^{12},r^{13}] = 
\sum_{i,j,\beta}  D_{ij} \phi_{\beta}  [x_i,e_{\beta}] \otimes
x_j \otimes e_{-\beta}  +
\sum_{k,l,\alpha} D_{kl} \phi_{\alpha} [e_{\alpha},x_k] \otimes
e_{-\alpha} \otimes x_l \\+
\sum_{\alpha,\beta} (-1)^{|\alpha||\beta|}
\phi_{\alpha} \phi_{\beta} [e_{\alpha},e_{\beta}] \otimes
e_{-\alpha} \otimes e_{-\beta}
\end{eqnarray*}
\begin{eqnarray*}
[r^{12},r^{23}] = 
\sum_{i,j,\beta}  D_{ij} \phi_{\beta}  x_i \otimes
[x_j,e_{\beta}] \otimes e_{-\beta}  +
\sum_{k,l,\alpha} D_{kl} \phi_{\alpha} e_{\alpha} \otimes
[e_{-\alpha},x_k] \otimes x_l \\+
\sum_{\alpha,\beta} 
\phi_{\alpha} \phi_{\beta} e_{\alpha} \otimes
[e_{-\alpha}, e_{\beta}] \otimes e_{-\beta}
\end{eqnarray*}
\begin{eqnarray*}
[r^{13},r^{23}] =
\sum_{i,j,\beta}  D_{ij} \phi_{\beta}  x_i \otimes e_{\beta}
\otimes [x_j, e_{-\beta}]  +
\sum_{k,l,\alpha} D_{kl} \phi_{\alpha} e_{\alpha} \otimes x_k
\otimes [e_{-\alpha}, x_l] \\+
\sum_{\alpha,\beta} (-1)^{|\alpha||\beta|}
\phi_{\alpha} \phi_{\beta} e_{\alpha} \otimes e_{\beta} \otimes
[e_{-\alpha}, e_{-\beta}]
\end{eqnarray*}

We can now easily see that the $\h \otimes \g_{\alpha} \otimes
\g_{-\alpha}$ part will be:
\begin{eqnarray*}
Alt_s(dr) \left|_{\h \otimes \g_{\alpha} \otimes \g_{-\alpha}}
\right. +
\sum_{\alpha,\beta} (-1)^{|\alpha||\beta|} \delta_{\alpha,-\beta}
\phi_{\alpha} \phi_{\beta} [e_{\alpha},e_{\beta}] \otimes
e_{-\alpha} \otimes e_{-\beta} + \\
\sum_{i,j,\beta}  D_{ij} \phi_{\beta}  x_i \otimes
[x_j,e_{\beta}] \otimes e_{-\beta}  +
\sum_{i,j,\beta}  D_{ij} \phi_{\beta}  x_i \otimes e_{\beta}
\otimes [x_j, e_{-\beta}]  
\end{eqnarray*}

We note that the last two sums cancel out, because:
$$ e_{\beta} \otimes [x_j, e_{-\beta}] = (-\beta)(x_j)
e_{\beta} \otimes e_{-\beta} = -(\beta(x_j)) e_{\beta}
\otimes e_{-\beta} = - [x_j, e_{\beta}] \otimes e_{\beta} $$ 
\noindent
and we have:
\begin{eqnarray*} 
&&Alt_s(dr) \left|_{\h \otimes \g_{\alpha} \otimes
\g_{-\alpha}} \right. + 
\sum_{\alpha,\beta} (-1)^{|\alpha||\beta|}
\delta_{\alpha,-\beta} \phi_{\alpha} \phi_{\beta}
[e_{\alpha},e_{\beta}] \otimes e_{-\alpha} \otimes e_{-\beta} \\
&=& 
\sum_{i,\alpha} \frac{\partial \phi_{\alpha}}{\partial x_i}
x_i \otimes e_{\alpha} \otimes e_{-\alpha} +
\sum_{\alpha} (-1)^{|\alpha|} 
\phi_{\alpha} \phi_{-\alpha} [e_{\alpha},e_{-\alpha}] \otimes
e_{-\alpha} \otimes e_{\alpha} \\
&=& \sum_{i,\alpha} 
\frac{\partial \phi_{\alpha}}{\partial x_i} x_i \otimes
e_{\alpha} \otimes e_{-\alpha} +
\sum_{\alpha} (-1)^{|\alpha|} 
\phi_{\alpha} \phi_{-\alpha} [e_{-\alpha},e_{\alpha}] \otimes
e_{\alpha} \otimes e_{-\alpha} \\
&=& \sum_{i,\alpha} 
\frac{\partial \phi_{\alpha}}{\partial x_i} x_i \otimes
e_{\alpha} \otimes e_{-\alpha} +
\sum_{\alpha} (-1)^{|\alpha|} 
\phi_{\alpha} \phi_{-\alpha} (-A_{\alpha}) h_{\alpha} \otimes
e_{\alpha} \otimes e_{-\alpha} \\
&=& \sum_{i,\alpha} 
\frac{\partial \phi_{\alpha}}{\partial x_i} x_i \otimes
e_{\alpha} \otimes e_{-\alpha} +
\sum_{\alpha} A_{\alpha}
\phi_{\alpha} \phi_{\alpha}  h_{\alpha} \otimes
e_{\alpha} \otimes e_{-\alpha} 
\end{eqnarray*}
\noindent
which we can rewrite as:
$$\sum_{\alpha \in \Delta}
\left( \sum_{i} 
\frac{\partial \phi_{\alpha}}{\partial x_i} x_i 
+ A_{\alpha} \phi_{\alpha}^2 h_{\alpha} 
\right) \otimes e_{\alpha} \otimes e_{-\alpha}
$$
\noindent
where $h_{\alpha} \in \h$ is defined by:
$$ [e_{\alpha},e_{-\alpha}] = (e_{\alpha},e_{-\alpha})h_{\alpha}
= A_{-\alpha} h_{\alpha}.$$
\noindent
(Recall that both $\{x_i\}$ and $\{h_{\alpha}|
\alpha \in \Gamma\}$ are bases for $\h$ and equivalently linear
coordinate systems for $\h^*$). For this term to vanish we must
have, for all ${\alpha \in \Delta}$:  $$ \sum_{i} 
\frac{\partial \phi_{\alpha}}{\partial x_i} x_i 
+  A_{\alpha} \phi_{\alpha}^2 h_{\alpha} = 0.$$
\noindent
We can rewrite this as: 
$$ d\phi_{\alpha} + A_{\alpha} \phi_{\alpha}^2 dh_{\alpha}
= 0. $$
\noindent
If we define $\mu_{\alpha}$ by: 
$$ \mu_{\alpha} = \left\{ 
\begin{array}{cl}
\sqrt{-1} & \textmd{if } \alpha \textmd{ is an odd positive root}
\\  1 & \textmd{otherwise}
\end{array}  \right.$$
\noindent
then $\mu_{\alpha}^2 = A_{\alpha}$, and the
equation we need to solve is: 
$$ d\phi_{\alpha} + \mu_{\alpha}^2 \phi_{\alpha}^2 dh_{\alpha}
= 0. $$
We assume $\phi_{\alpha} \neq 0$ and let $u_{\alpha}=
\mu_{\alpha}\phi_{\alpha}$. Separating variables to
integrate we obtain:  
$$ \frac{1}{\mu_{\alpha}} 
\int \frac{-du_{\alpha}}{u_{\alpha}^2} = \int dh_{\alpha} 
\; \; \; \Rightarrow \; \; \;  
u_{\alpha} = 
\mu_{\alpha} \phi_{\alpha} = 
\frac{1}{\mu_{\alpha} h_{\alpha} + C} 
$$  
\noindent
and we get: 
$$ \phi_{\alpha} = \frac{ A_{\alpha}}{h_{\alpha} -
\nu_{\alpha}}$$
\noindent
for some ${\nu_{\alpha} \in \C}$. Here $h_{\alpha}$ is viewed as
a linear function on $\h^*$ via $h_{\alpha}(\lambda) = {(\alpha,
\lambda)}$.

Finally we look at the $\g_{\alpha} \otimes \g_{\beta} \otimes
\g_{-\alpha-\beta}$ part of Equation \ref{SDYBequation}. There
is no contribution from the dynamical part; the only terms we
need to look at are: 
\begin{eqnarray*}
\sum_{\alpha,\beta} (-1)^{|\alpha||\beta|}
\phi_{\alpha} \phi_{\beta} [e_{\alpha},e_{\beta}] \otimes
e_{-\alpha} \otimes e_{-\beta} +
\sum_{\alpha,\beta} 
\phi_{\alpha} \phi_{\beta} e_{\alpha} \otimes
[e_{-\alpha}, e_{\beta}] \otimes e_{-\beta}
\\ +
\sum_{\alpha,\beta} (-1)^{|\alpha||\beta|}
\phi_{\alpha} \phi_{\beta} e_{\alpha} \otimes e_{\beta} \otimes
[e_{-\alpha}, e_{-\beta}]
\end{eqnarray*}
\noindent
Denote by $C$ the matrix of coefficients determined by:
$$ [e_{\alpha}, e_{\beta}] = C_{\alpha,\beta}^{\gamma}
e_{\gamma}.$$
\noindent
Then we can rewrite the terms we are interested in as: 
\begin{eqnarray*}
&\sum_{\alpha,\beta,\gamma} (-1)^{|\alpha||\beta|}
\phi_{\alpha} \phi_{\beta} C_{\alpha,\beta}^{\gamma} e_{\gamma}
\otimes e_{-\alpha} \otimes e_{-\beta} \\ &\qquad +
\sum_{\alpha,\beta,\gamma} 
\phi_{\alpha} \phi_{\beta} C_{-\alpha,\beta}^{\gamma} e_{\alpha}
\otimes e_{\gamma} \otimes e_{-\beta} \\ &
\qquad \qquad +
\sum_{\alpha,\beta,\gamma} (-1)^{|\alpha||\beta|}
\phi_{\alpha} \phi_{\beta} C_{-\alpha,-\beta}^{\gamma} e_{\alpha}
\otimes e_{\beta} \otimes e_{\gamma}  \\ \\ &=
\sum_{\alpha,\beta} (-1)^{|\alpha||\beta|}
\phi_{\alpha} \phi_{\beta} C_{\alpha,\beta}^{\alpha+\beta}
e_{\alpha+\beta} \otimes e_{-\alpha} \otimes e_{-\beta} \\ &\qquad +
\sum_{\alpha,\beta} 
\phi_{\alpha} \phi_{\beta} C_{-\alpha,\beta}^{-\alpha+\beta}
e_{\alpha} \otimes e_{-\alpha+\beta} \otimes e_{-\beta} \\ &
\qquad \qquad +
\sum_{\alpha,\beta} (-1)^{|\alpha||\beta|}
\phi_{\alpha} \phi_{\beta} C_{-\alpha,-\beta}^{-\alpha-\beta}
e_{\alpha} \otimes e_{\beta} \otimes e_{-\alpha-\beta} \\ \\ &= 
\sum_{\alpha,\beta} (-1)^{|\alpha||\beta|}
\phi_{\alpha} \phi_{\beta} C_{\alpha,\beta}^{\alpha+\beta}
e_{\alpha+\beta} \otimes e_{-\alpha} \otimes e_{-\beta} \\ &\qquad +
\sum_{\alpha,\beta} 
\phi_{\alpha} \phi_{\beta} C_{-\alpha,\beta}^{-\alpha+\beta}
e_{\alpha} \otimes e_{-\alpha+\beta} \otimes e_{-\beta} \\ &
\qquad \qquad + 
\sum_{\alpha,\beta} (-1)^{|\alpha||\beta|}
\phi_{\alpha} \phi_{\beta} C_{-\alpha,-\beta}^{-\alpha-\beta}
e_{\alpha} \otimes e_{\beta} \otimes e_{-\alpha-\beta} 
\end{eqnarray*}
\noindent
We want the coefficient in front of the term ${e_{\alpha}
\otimes e_{\beta} \otimes e_{-\alpha-\beta}}$ to vanish:
\begin{eqnarray*}
(-1)^{|\beta||\alpha+\beta|} \phi_{-\beta} \phi_{\alpha+\beta}
C_{-\beta,\alpha+\beta}^{\alpha} 
+ \phi_{\alpha} \phi_{\alpha+\beta}
C_{-\alpha,\alpha+\beta}^{\beta} &
\\ \\ + (-1)^{|\alpha||\beta|} \phi_{\alpha} \phi_{\beta}
C_{-\alpha,-\beta}^{-\alpha-\beta} &
= \\ \\  
-(-1)^{|\alpha||\beta|} \phi_{\beta} \phi_{\alpha+\beta}
C_{-\beta,\alpha+\beta}^{\alpha} 
+ \phi_{\alpha} \phi_{\alpha+\beta}
C_{-\alpha,\alpha+\beta}^{\beta} &
\\ \\ + (-1)^{|\alpha||\beta|} \phi_{\alpha} \phi_{\beta}
C_{-\alpha,-\beta}^{-\alpha-\beta} &
=0
\end{eqnarray*}
\noindent
We now compute these constants: 
\begin{eqnarray*}
C_{-\beta,\alpha+\beta}^{\alpha} &=&
([e_{-\beta},e_{\alpha+\beta}], e_{\alpha}^*) 
= ([e_{-\beta},e_{\alpha+\beta}],A_{-\alpha}e_{-\alpha}) \\ \\
&=&
\frac{A_{-\alpha}}{C_{\alpha,\beta}^{\alpha+\beta}}
([e_{-\beta},[e_{\alpha},e_{\beta}]],e_{-\alpha}) \\
&=& 
\frac{-(-1)^{|\alpha||\beta|}A_{-\alpha}}
{C_{\alpha,\beta}^{\alpha+\beta}}
([e_{-\beta},[e_{\beta},e_{\alpha}]],e_{-\alpha}) \\
&=& 
\frac{-(-1)^{|\alpha||\beta|}A_{-\alpha}}
{C_{\alpha,\beta}^{\alpha+\beta}}
([[e_{-\beta},e_{\beta}],e_{\alpha}] +
(-1)^{|\beta|}[e_{\beta}.[e_{-\beta},e_{\alpha}]],e_{-\alpha}) 
\\
&=& 
\frac{-(-1)^{|\alpha||\beta|}A_{-\alpha}}
{C_{\alpha,\beta}^{\alpha+\beta}}
([[e_{-\beta},e_{\beta}],e_{\alpha}],e_{-\alpha}) \\
&=& 
\frac{(-1)^{|\alpha||\beta|}(-1)^{|\beta|}A_{-\alpha}}
{C_{\alpha,\beta}^{\alpha+\beta}}
([[e_{\beta},e_{-\beta}],e_{\alpha}],e_{-\alpha}) 
\\
&=& 
\frac{(-1)^{|\alpha||\beta|}(-1)^{|\beta|}A_{-\alpha}A_{-\beta}}
{C_{\alpha,\beta}^{\alpha+\beta}}
(h_{\beta},[e_{\alpha},e_{-\alpha}]) \\
&=& 
\frac{(-1)^{|\alpha||\beta|}(-1)^{|\beta|}
A_{-\beta}}{C_{\alpha,\beta}^{\alpha+\beta}} (h_{\alpha},h_{\beta}) 
\end{eqnarray*} 
\noindent
where we use the Jacobi identity and $(A_{\alpha})^2 = 1$.
Similar computations yield: 
\begin{eqnarray*}
C_{-\alpha,\alpha+\beta}^{\beta} 
&=& 
\frac{-(-1)^{|\alpha|}A_{-\alpha}}
{C_{\alpha,\beta}^{\alpha+\beta}} (h_{\alpha},h_{\beta}) \\
C_{-\alpha,-\beta}^{-\alpha-\beta} 
&=& 
\frac{(-1)^{|\alpha||\beta|}(-1)^{|\alpha|+|\beta|}A_{\alpha+\beta}
A_{-\alpha}A_{-\beta}}
{C_{\alpha,\beta}^{\alpha+\beta}} (h_{\alpha},h_{\beta}) 
\end{eqnarray*}
\noindent
Recall that if $\alpha+\beta$ is a root,
$(h_{\alpha},h_{\beta}) \neq 0$ and
$C_{\alpha,\beta}^{\alpha+\beta} \neq 0$. Then we can rewrite
the coefficient of ${e_{\alpha} \otimes e_{\beta} \otimes
e_{-\alpha-\beta}}$: 
\begin{eqnarray*}
0&=& -(-1)^{|\alpha||\beta|} \phi_{\beta} \phi_{\alpha+\beta}
C_{-\beta,\alpha+\beta}^{\alpha} 
+ \phi_{\alpha} \phi_{\alpha+\beta}
C_{-\alpha,\alpha+\beta}^{\beta} \\ && \qquad \qquad
+ (-1)^{|\alpha||\beta|} \phi_{\alpha} \phi_{\beta}
C_{-\alpha,-\beta}^{-\alpha-\beta} \\ \\
0&=& -(-1)^{|\alpha||\beta|} \phi_{\beta} \phi_{\alpha+\beta}
\frac{(-1)^{|\alpha||\beta|}(-1)^{|\beta|}A_{-\beta}}
{C_{\alpha,\beta}^{\alpha+\beta}} (h_{\alpha},h_{\beta}) \\ 
&& + \;\;\; \phi_{\alpha} \phi_{\alpha+\beta}
\frac{-(-1)^{|\alpha|}A_{-\alpha}}
{C_{\alpha,\beta}^{\alpha+\beta}} (h_{\alpha},h_{\beta}) \\ 
&& + \;\;\; (-1)^{|\alpha||\beta|} \phi_{\alpha} \phi_{\beta}
\frac{(-1)^{|\alpha||\beta|}(-1)^{|\alpha|+|\beta|}A_{\alpha+\beta}
A_{-\alpha}A_{-\beta}}
{C_{\alpha,\beta}^{\alpha+\beta}} (h_{\alpha},h_{\beta}) \\ \\
0&=& -(-1)^{|\alpha||\beta|} \phi_{\beta} \phi_{\alpha+\beta}
(-1)^{|\alpha||\beta|}(-1)^{|\beta|}A_{-\alpha}
 + \phi_{\alpha} \phi_{\alpha+\beta}
(-(-1)^{|\alpha|})A_{-\beta} \\
&& + \;\;\; (-1)^{|\alpha||\beta|} \phi_{\alpha} \phi_{\beta}
(-1)^{|\alpha||\beta|}(-1)^{|\alpha|+|\beta|}A_{\alpha+\beta} \\ \\
0&=& -(-1)^{|\beta|} \phi_{\beta} \phi_{\alpha+\beta} A_{-\alpha}
- (-1)^{|\alpha|} \phi_{\alpha} \phi_{\alpha+\beta}
A_{-\beta} \\ && \qquad \qquad
+ (-1)^{|\alpha|+|\beta|} \phi_{\alpha} \phi_{\beta}
A_{\alpha+\beta} 
\end{eqnarray*}
\noindent
which can be rewritten as:
\begin{eqnarray*}
A_{\alpha+\beta}\phi_{\alpha}\phi_{\beta} &=& 
A_{\alpha}\phi_{\beta}\phi_{\alpha+\beta} +
A_{\beta}\phi_{\alpha}\phi_{\alpha+\beta}\\
&=& \left(A_{\alpha}\phi_{\beta} + A_{\beta}\phi_{\alpha} \right)
\phi_{\alpha+\beta}
\end{eqnarray*}
\noindent
where we use $A_{\alpha}= (-1)^{|\alpha|}A_{-\alpha}$ and
$(A_{\alpha})^2 = 1$.

Define the set 
$$X = \{\alpha \in \Delta | \phi_{\alpha} \neq 0 \}$$
\noindent
Then it is easy to see that $X$ is closed under addition and
changing signs. In other words, if ${\alpha, \beta \in X,}$ then 
so are $\alpha+\beta$ and $-\alpha$. These
follow directly from the above equation relating
$\phi_{\alpha}$ and $\phi_{\beta}$ to $\phi_{\alpha+\beta}$, and
the unitarity property (i.e. $r = -T_s(r)$).

Next assume $\alpha$ and $\beta$ are positive and in $X$. Then
the above calculations yield:
\begin{eqnarray*}
A_{\alpha+\beta}\phi_{\alpha}\phi_{\beta} 
=& \left(A_{\alpha}\phi_{\beta} + A_{\beta}\phi_{\alpha} \right)
\phi_{\alpha+\beta}\\ \\
(-1)^{|\alpha|+|\beta|}\phi_{\alpha}\phi_{\beta} 
=& \left((-1)^{|\alpha|}\phi_{\beta} +
(-1)^{|\beta|}\phi_{\alpha} \right) \phi_{\alpha+\beta}\\ \\
(-1)^{|\alpha|+|\beta|}
\left(\frac{(-1)^{|\alpha|}}{h_{\alpha}-\nu_{\alpha}} \right)
\left(\frac{(-1)^{|\beta|}}{h_{\beta}-\nu_{\beta}}\right)
=& \left(
\frac{(-1)^{|\alpha+\beta|}}{h_{\beta}-\nu_{\beta}}
+ 
\frac{(-1)^{|\alpha+\beta|}}{h_{\alpha}-\nu_{\alpha}}
 \right)
\frac{(-1)^{|\alpha+\beta|}}{h_{\alpha+\beta}-\nu_{\alpha+\beta}}
\\ \\
\left(\frac{1}{h_{\alpha}-\nu_{\alpha}} \right)
\left(\frac{1}{h_{\beta}-\nu_{\beta}}\right)
=& \left(
\frac{1}{h_{\beta}-\nu_{\beta}}
+ 
\frac{1}{h_{\alpha}-\nu_{\alpha}}
 \right)
\frac{1}{h_{\alpha+\beta}-\nu_{\alpha+\beta}}
\end{eqnarray*}
\noindent
Since $h_{\alpha+\beta} = \h_{\alpha} + h_{\beta}$ we must have:
$$\nu_{\alpha+\beta} = \nu_{\alpha} +\nu_{\beta}$$
\noindent
Also recall that $\phi_{-\alpha} = -(-1)^{|\alpha|}\phi_{\alpha}$
and so:
$$\frac{A_{-\alpha}}{h_{-\alpha} - \nu_{-\alpha}}
=  -(-1)^{|\alpha|}\frac{A_{\alpha}}{h_{\alpha}-\nu_{\alpha}} $$ 
\noindent 
which implies that $\nu_{-\alpha} = -\nu_{\alpha}$.
Therefore we can conclude that there is some $\nu \in \h^*$ such
that $\nu_{\alpha} = (\alpha,\nu)$ for all $\alpha \in X$. 
This completes the proof of the theorem.
$\blacksquare$

\subsection{Proof of Theorem \ref{0weighttheorem}}
\label{Proofof0Thm}

We start with fixing a basis for the non-Cartan part of $\g$ in
the same manner as above. In other words for any positive root
$\alpha$ we fix $e_{\alpha} \in \g_{\alpha}$ and pick
${e_{-\alpha} \in \g_{-\alpha}}$ dual to $e_{\alpha}$ i.e. 
$$ (e_{\alpha}, e_{-\alpha}) = 1 \textmd{ for all } \alpha \in
\Delta^+.$$
\noindent
We again need the following notation:
$$A_{\alpha} = \left\{ \begin{array}{cl} 
(-1)^{|\alpha|} & \textmd{if } \alpha  \textmd{ is positive} \\
1 & \textmd{if } \alpha \textmd{ is negative}
\end{array} \right. $$
\noindent
and we note once again that $A_{-\alpha} =
(-1)^{|\alpha|}A_{\alpha}$.

Let $r : \h^* \rightarrow \g \otimes \g$ be a meromorphic map,
$\Omega \in \g \otimes \g$ the Casimir element,
and $\epsilon$ a nonzero complex number. Introduce a second
meromorphic function ${s : \h^* \rightarrow \g \otimes \g}$ by:
$$ s(\lambda) = r(\lambda) - \frac{\epsilon}{2}\Omega  \;
\;\;\;\; \textmd{ for all } \lambda \; \in \; \h^*.$$

We will now prove the following technical lemma:

\vspace{0.1in}
\begin{lemma}
The map $r$ is a super dynamical r-matrix with zero
weight and coupling constant $\epsilon$ if and only if $s$
satisfies the zero weight condition:
$$ [h\otimes 1 + 1 \otimes h, s(\lambda)] = 0 \textmd{ for all }
h \in \h, \lambda \in \h^*,$$
\noindent
the unitarity condition:
$$ s(\lambda) + T_s(s)(\lambda) = 0,$$
\noindent
and the following modified version of the dynamical Yang-Baxter
equation: 
\begin{equation}
\label{MDYBequation}
Alt_s(ds) + [[s,s]] + \frac{\epsilon^2}{4}[[\Omega,\Omega]] = 0,
\end{equation}
\noindent
where, for any $2-$tensor $r$:
$$ [[r,r]] =  [r^{12},r^{13}] + [r^{12},r^{23}] +
[r^{13},r^{23}].$$
\end{lemma}

\vspace{0.1in}
\begin{remark}
This is exactly Lemma $3.9$ of \cite{EV}. Its proof consists
mainly of a direct calculation,  but the computations are
significantly more involved in the super case. In any case,
the proof will be given here for completeness.
\end{remark}

\vspace{0.1in}
\noindent
\textbf{Proof: }
The zero weight condition and the unitarity condition on $s$ 
imply that $s$ has to be of the form:
$$ s(\lambda) = \sum_{i,j} D_{ij}(\lambda) h_i \otimes h_j +
\sum_{\alpha \in \Delta} \phi_{\alpha} e_{\alpha} \otimes
e_{-\alpha}, $$
\noindent 
where $\{h_i\}$ is a basis for $\h$, and $D_{ij}, \phi_{\alpha}$
are suitable scalar meromorphic functions such that:
$$ D_{ij}(\lambda) = - D_{ji}(\lambda)  
\; \; \; \; \;
\textmd{and} 
\; \; \; \; \;
\phi_{-\alpha} = -(-1)^{|\alpha|} \phi_{\alpha}$$
\noindent
One can easily check that there can be no terms mixing the
Cartan part with the non-Cartan part, and the conditions on the
$D_{ij}$ and the $\phi_{\alpha}$ follow from unitarity.

We have:
$$ [[r,r]] = [[s + \frac{\epsilon}{2}\;\Omega, s +
\frac{\epsilon}{2}\;\Omega\;]] $$
\noindent
so it suffices to prove that $[[r,r]] -  [[s,s]] -
\frac{\epsilon^2}{4}[[\;\Omega, \Omega\;]]$ is zero. In other
words we need to show that  
$$\frac{\epsilon}{2}  \left(
[s^{12},\Omega^{13}] + [\Omega^{12},s^{13}] +
[s^{12},\Omega^{23}] + [\Omega^{12},s^{23}] +
[s^{13},\Omega^{23}] + [\Omega^{13},s^{23}] \right) $$
\noindent
vanishes.

If $\{h^i\}$ is the dual basis in $\h$ to $\{h_i\}$, then we
can write the Casimir element as:

$$ \Omega = \sum_i h_i \otimes h^i + 
\sum_{\alpha \in \Delta} A_{\alpha} e_{\alpha} \otimes
e_{-\alpha} $$

\noindent
and we have:

\begin{eqnarray*} 
\left[ s^{12},\Omega^{13}\right] &=& 
\sum_{i,j,k}D_{ij} [h_i, h_k] \otimes h_j \otimes h^k +
\sum_{\alpha,k} \phi_{\alpha}[e_{\alpha},h_k] \otimes e_{-\alpha}
\otimes h^k \\
&+&
\sum_{i,j,\beta} D_{ij}A_{\beta} [h_i,e_{\beta}] \otimes h_j
\otimes e_{-\beta} \\ && \qquad +
\sum_{\alpha,\beta} (-1)^{|\alpha||\beta|}
\phi_{\alpha}A_{\beta}[e_{\alpha},e_{\beta}]\otimes
e_{-\alpha} \otimes e_{-\beta}  
\\ 
\left[ \Omega^{12},s^{13}\right] &=& 
\sum_{i,j,k} D_{jk}[h_i,h_j]\otimes h^i \otimes h_k +
\sum_{i,\beta} \phi_{\beta}[h_i,e_{\beta}] \otimes h^i \otimes
e_{-\beta} \\ 
&+&
\sum_{\alpha,j,k} A_{\alpha}D_{jk} [e_{\alpha},h_j] \otimes
e_{-\alpha} \otimes h_k \\ && \qquad +
\sum_{\alpha,\beta} (-1)^{|\alpha||\beta|}
A_{\alpha}\phi_{\beta} [e_{\alpha},e_{\beta}]
\otimes e_{-\alpha} \otimes e_{-\beta} 
\\ 
\left[ s^{12},\Omega^{23}\right] &=&
\sum_{i,j,k}D_{ij} h_i \otimes [h_j,h_k] \otimes h^k +
\sum_{\alpha,k} \phi_{\alpha}e_{\alpha} \otimes [e_{-\alpha},h_k]
\otimes h^k \\
&+&
\sum_{i,j,\beta} D_{ij}A_{\beta} h_i \otimes [h_j,e_{\beta}]
\otimes e_{-\beta} \\ && \qquad +
\sum_{\alpha,\beta}\phi_{\alpha}A_{\beta} e_{\alpha} \otimes
[e_{-\alpha},e_{\beta}] \otimes e_{-\beta}  
\\ 
\left[ \Omega^{12},s^{23}\right] &=&
\sum_{i,j,k} D_{jk}h_i \otimes [h^i,h_j] \otimes h_k +
\sum_{i,\beta} \phi_{\beta} h_i \otimes [h^i,e_{\beta}] \otimes
e_{-\beta} \\ 
&+&
\sum_{\alpha,j,k} A_{\alpha}D_{jk} e_{\alpha} \otimes
[e_{-\alpha}, h_j] \otimes h_k \\ && \qquad +
\sum_{\alpha,\beta} A_{\alpha}\phi_{\beta} e_{\alpha}
\otimes [e_{-\alpha},e_{\beta}] \otimes e_{-\beta} 
\\ 
\left[ s^{13},\Omega^{23}\right] &=& 
\sum_{i,j,k}D_{ij} h_i \otimes h_k \otimes [h_j, h^k] +
\sum_{\alpha,k} \phi_{\alpha}e_{\alpha} \otimes h_k \otimes
[e_{-\alpha}, h^k] \\
&+&
\sum_{i,j,\beta} D_{ij}A_{\beta} h_i \otimes e_{\beta} \otimes
[h_j, e_{-\beta}] \\ && \qquad +
\sum_{\alpha,\beta} (-1)^{|\alpha||\beta|}
\phi_{\alpha}A_{\beta}
e_{\alpha} \otimes e_{\beta} \otimes [e_{-\alpha}, e_{-\beta}]  
\\ 
\left[ \Omega^{13},s^{23}\right] &=& 
\sum_{i,j,k} D_{jk}h_i \otimes h_j \otimes [h^i, h_k] +
\sum_{i,\beta} \phi_{\beta} h_i \otimes e_{\beta} \otimes [h^i,
e_{-\beta}] \\  
&+&
\sum_{\alpha,j,k} A_{\alpha}D_{jk} e_{\alpha} \otimes h_j
\otimes [e_{-\alpha}, h_k] \\ && \qquad +
\sum_{\alpha,\beta} (-1)^{|\alpha||\beta|}
A_{\alpha}\phi_{\beta} e_{\alpha}
\otimes e_{\beta} \otimes [e_{-\alpha}, e_{-\beta}] 
\end{eqnarray*}

Clearly the first sum in each of these six terms is zero, as the
$h_i$ are all in the Cartan. It is also easy to see that the
sum of all the third sums vanishes as well: 
\begin{eqnarray*} 
\sum_{i,j,\beta} D_{ij}A_{\beta} [h_i,e_{\beta}] \otimes h_j
\otimes e_{-\beta} &=& 
\sum_{i,j,\beta} D_{ij}A_{\beta} \beta(h_i) e_{\beta} \otimes
h_j \otimes e_{-\beta} \\
&=& 
-\sum_{i,j,\beta} (-D_{ij})A_{\beta} e_{\beta} \otimes
h_j \otimes [e_{-\beta},h_i] \\
&=&
-\sum_{\alpha,j,k} D_{j,k}A_{\alpha} e_{\alpha} \otimes
h_j \otimes [e_{-\alpha},h_k]
\end{eqnarray*}
\noindent
where we use $D_{ij} = -D_{ji}$,
\begin{eqnarray*}
\sum_{\alpha,j,k} A_{\alpha}D_{jk} [e_{\alpha},h_j] \otimes
e_{-\alpha} \otimes h_k  &=& 
-\sum_{\alpha,j,k} A_{\alpha}D_{jk} \alpha(h_j) e_{\alpha}
\otimes e_{-\alpha} \otimes h_k \\
&=&
-\sum_{\alpha,j,k} A_{\alpha}D_{jk} e_{\alpha} \otimes
[e_{-\alpha}, h_j] \otimes h_k, 
\end{eqnarray*}
\noindent
and:
\begin{eqnarray*}
\sum_{i,j,\beta} D_{ij}A_{\beta} h_i \otimes [h_j,e_{\beta}]
\otimes e_{-\beta} &=&
\sum_{i,j,\beta} D_{ij}A_{\beta} \beta(h_j) h_i \otimes e_{\beta}
\otimes e_{-\beta} \\
&=&
-\sum_{i,j,\beta} D_{ij}A_{\beta} h_i \otimes e_{\beta} \otimes
[h_j, e_{-\beta}].
\end{eqnarray*}
\noindent
We can also see that:
\begin{eqnarray*} 
\sum_{\alpha,k} \phi_{\alpha}[e_{\alpha},h_k] \otimes e_{-\alpha}
\otimes h^k  &=& 
\sum_{\alpha,k} \phi_{\alpha} (-\alpha(h_k) e_{\alpha}) \otimes
e_{-\alpha} \otimes h^k \\
&=& 
-\sum_{\alpha,k} \phi_{\alpha} e_{\alpha} \otimes [e_{-\alpha},
h_k] \otimes h^k,
\end{eqnarray*} 
\noindent
and:
\begin{eqnarray*} 
\sum_{i,\beta} \phi_{\beta} h_i \otimes [h^i,e_{\beta}] \otimes
e_{-\beta} &=& 
\sum_{i,\beta} \phi_{\beta} \beta(h^i) h_i \otimes e_{\beta}
\otimes e_{-\beta} \\
&=& 
-\sum_{i,\beta} \phi_{\beta} h_i \otimes e_{\beta} \otimes [h^i,
e_{-\beta}] 
\end{eqnarray*}
\noindent
but: 
\begin{eqnarray*}
\sum_{i,\beta} \phi_{\beta}[h_i,e_{\beta}] \otimes h^i \otimes
e_{-\beta} &=& 
\sum_{i,\beta} \phi_{\beta} \beta(h_i) e_{\beta} \otimes h^i
\otimes e_{-\beta} \\ 
&=&
\sum_{\beta} \phi_{\beta} e_{\beta} \otimes (\sum_i
\beta(h_i)h^i) \otimes e_{-\beta} \\
&=&
\sum_{\beta} \phi_{\beta} e_{\beta} \otimes h_{\beta} \otimes
e_{-\beta}\\ 
\sum_{\alpha,k} \phi_{\alpha}e_{\alpha} \otimes h_k
\otimes [e_{-\alpha}, h^k] &=&
\sum_{\alpha,k} \phi_{\alpha} \alpha(h^k) e_{\alpha} \otimes h_k
\otimes e_{-\alpha} \\
&=&
\sum_{\alpha} \phi_{\alpha} e_{\alpha} \otimes
(\sum_k \alpha(h^k) h_k) \otimes e_{-\alpha} \\
&=& 
\sum_{\alpha} \phi_{\alpha} e_{\alpha} \otimes h_{\alpha}
\otimes e_{-\alpha}
\end{eqnarray*}
\noindent
where $h_{\alpha} \in \h$ is
defined as usual by
${A_{-\alpha}h_{\alpha} = (e_{\alpha},e_{-\alpha})h_{\alpha} =
[e_{\alpha},e_{-\alpha}].}$

Next we look at the six remaining terms, the fourth sum in each
term above. We note that: 
\begin{eqnarray*} 
\sum_{\alpha,\beta} (-1)^{|\alpha||\beta|} 
\phi_{\alpha}A_{\beta}[e_{\alpha},e_{\beta}] \otimes
e_{-\alpha} &\otimes& e_{-\beta} \delta_{\alpha,-\beta} \\
&=&
\sum_{\alpha} (-1)^{|\alpha|} \phi_{\alpha} A_{-\alpha}
[e_{\alpha},e_{-\alpha}] \otimes
e_{-\alpha} \otimes e_{\alpha}
\\
\sum_{\alpha,\beta} (-1)^{|\alpha||\beta|}
A_{\alpha}\phi_{\beta} [e_{\alpha},e_{\beta}]
\otimes e_{-\alpha} &\otimes& e_{-\beta} \delta_{\alpha,-\beta} \\
&=&
\sum_{\alpha} (-1)^{|\alpha|}
A_{\alpha}\phi_{-\alpha} [e_{\alpha},e_{-\alpha}]
\otimes e_{-\alpha} \otimes e_{\alpha}  
\\
\sum_{\alpha,\beta}\phi_{\alpha}A_{\beta} e_{\alpha} \otimes
[e_{-\alpha},e_{\beta}] &\otimes& e_{-\beta} 
\delta_{\alpha,\beta} \\
&=&
\sum_{\alpha}\phi_{\alpha}A_{\alpha} e_{\alpha} \otimes
[e_{-\alpha},e_{\alpha}] \otimes e_{-\alpha}  
\\
\sum_{\alpha,\beta} A_{\alpha}\phi_{\beta} e_{\alpha}
\otimes [e_{-\alpha},e_{\beta}] &\otimes& e_{-\beta}
\delta_{\alpha,\beta} \\
&=&
\sum_{\alpha} A_{\alpha}\phi_{\alpha} e_{\alpha}
\otimes [e_{-\alpha},e_{\alpha}] \otimes e_{-\alpha}
\\
\sum_{\alpha,\beta} (-1)^{|\alpha||\beta|}
\phi_{\alpha}A_{\beta} e_{\alpha} \otimes e_{\beta} &\otimes&
[e_{-\alpha}, e_{-\beta}]   \delta_{\alpha,-\beta} \\
&=&
\sum_{\alpha} (-1)^{|\alpha|}
\phi_{\alpha}A_{-\alpha} e_{\alpha} \otimes e_{-\alpha} \otimes
[e_{-\alpha}, e_{\alpha}] 
\\ 
\sum_{\alpha,\beta} (-1)^{|\alpha||\beta|}
A_{\alpha}\phi_{\beta} e_{\alpha} \otimes e_{\beta} &\otimes&
[e_{-\alpha}, e_{-\beta}]  \delta_{\alpha,-\beta} \\
&=&
\sum_{\alpha} (-1)^{|\alpha|}
A_{\alpha}\phi_{-\alpha} e_{\alpha} \otimes e_{-\alpha} \otimes
[e_{-\alpha}, e_{\alpha}]. 
\end{eqnarray*}
\noindent
We can see that the first two and the last two of these cancel
out one another because:
$$ \phi_{\alpha}A_{-\alpha} = (-1)^{|\alpha|}
\phi_{\alpha}A_{\alpha} \;\;\;\;\;\;
A_{\alpha}\phi_{-\alpha} =  -(-1)^{|\alpha|}
A_{\alpha}\phi_{\alpha}$$
\noindent
and the center two add up to give:
\begin{eqnarray*} 
2\sum_{\alpha}\phi_{\alpha}A_{\alpha} e_{\alpha} \otimes
[e_{-\alpha},e_{\alpha}] &\otimes& e_{-\alpha} \\
&=&
-2\sum_{\alpha}\phi_{\alpha}(-1)^{|\alpha|}A_{\alpha} e_{\alpha}
\otimes [e_{\alpha},e_{-\alpha}] \otimes e_{-\alpha}   
\\
&=&
-2\sum_{\alpha}\phi_{\alpha}(A_{-\alpha})^2 e_{\alpha}
\otimes h_{\alpha} \otimes e_{-\alpha}  
\\
&=&
-2\sum_{\alpha}\phi_{\alpha} e_{\alpha} \otimes h_{\alpha}
\otimes e_{-\alpha}   
\end{eqnarray*}
\noindent
because $A_{\alpha}^2 = 1$. This cancels the terms remaining
from the second sums:  
$$\sum_{\beta} \phi_{\beta} e_{\beta} \otimes h_{\beta} \otimes
e_{-\beta} +
\sum_{\alpha} \phi_{\alpha} e_{\alpha} \otimes h_{\alpha} \otimes
e_{-\alpha}
=
2\sum_{\beta} \phi_{\beta} e_{\beta} \otimes h_{\beta}
\otimes e_{-\beta}$$

Finally the only terms remaining that we need to check are the
terms of the form ${e_{\alpha}\otimes e_{\beta} \otimes
e_{\gamma}}$ which all come from the fourth sums. If $C$ is
again the matrix of coefficients determined by:
$$ [e_{\alpha},e_{\beta}] =
\sum_{\gamma} C_{\alpha,\beta}^{\gamma} e_{\gamma} =
C_{\alpha,\beta}^{\alpha+\beta} e_{\alpha+\beta}, $$ 
\noindent
then we can write the last six terms as follows:
\begin{eqnarray*} 
\sum_{\alpha,\beta} (-1)^{|\alpha||\beta|} 
\phi_{\alpha}A_{\beta}[e_{\alpha},e_{\beta}] &\otimes&
e_{-\alpha} \otimes e_{-\beta}  \\
&=&
\sum_{\alpha,\beta} (-1)^{|\alpha||\beta|} 
\phi_{\alpha}A_{\beta}
C_{\alpha,\beta}^{\alpha+\beta}e_{\alpha+\beta}\otimes
e_{-\alpha} \otimes e_{-\beta}   
\\
\sum_{\alpha,\beta} (-1)^{|\alpha||\beta|}
A_{\alpha}\phi_{\beta} [e_{\alpha},e_{\beta}]
&\otimes& e_{-\alpha} \otimes e_{-\beta} \\
&=&
\sum_{\alpha,\beta} (-1)^{|\alpha||\beta|}
A_{\alpha}\phi_{\beta}
C_{\alpha,\beta}^{\alpha+\beta}e_{\alpha+\beta} \otimes
e_{-\alpha} \otimes e_{-\beta}   
\\
\sum_{\alpha,\beta}\phi_{\alpha}A_{\beta} e_{\alpha}
&\otimes& [e_{-\alpha},e_{\beta}] \otimes e_{-\beta}  \\
&=&
\sum_{\alpha,\beta}\phi_{\alpha}A_{\beta} 
C_{-\alpha,\beta}^{-\alpha+\beta}
e_{\alpha} \otimes e_{-\alpha+\beta} \otimes e_{-\beta}  
\\
\sum_{\alpha,\beta} A_{\alpha}\phi_{\beta} e_{\alpha} &\otimes&
[e_{-\alpha},e_{\beta}] \otimes e_{-\beta}  \\
&=&
\sum_{\alpha,\beta} A_{\alpha}\phi_{\beta} 
C_{-\alpha,\beta}^{-\alpha+\beta} e_{\alpha}
\otimes e_{-\alpha+\beta} \otimes e_{-\beta} 
\\
\sum_{\alpha,\beta} (-1)^{|\alpha||\beta|}
\phi_{\alpha}A_{\beta}
e_{\alpha} &\otimes& e_{\beta} \otimes [e_{-\alpha}, e_{-\beta}]  \\
&=&
\sum_{\alpha,\beta} (-1)^{|\alpha||\beta|}
\phi_{\alpha}A_{\beta} C_{-\alpha,-\beta}^{-\alpha-\beta}
e_{\alpha} \otimes e_{\beta} \otimes e_{-\alpha-\beta}  
\\ 
\sum_{\alpha,\beta} (-1)^{|\alpha||\beta|}
A_{\alpha}\phi_{\beta} e_{\alpha}
&\otimes& e_{\beta} \otimes [e_{-\alpha}, e_{-\beta}] \\
&=&
\sum_{\alpha,\beta} (-1)^{|\alpha||\beta|}
A_{\alpha}\phi_{\beta} C_{-\alpha,-\beta}^{-\alpha-\beta}
e_{\alpha} \otimes e_{\beta} \otimes e_{-\alpha-\beta} 
\end{eqnarray*}
\noindent
We write down the contribution of each of these six sums to the
coefficient in front of ${e_{\alpha}\otimes e_{\beta} \otimes
e_{-\alpha-\beta}}$:
\begin{eqnarray*} 
\sum_{\alpha,\beta} (-1)^{|\alpha||\beta|} 
\phi_{\alpha}A_{\beta}
C_{\alpha,\beta}^{\alpha+\beta}e_{\alpha+\beta} &\otimes&
e_{-\alpha} \otimes e_{-\beta}   \\
&:&
(-1)^{|\beta||\alpha+\beta|} \phi_{-\beta}A_{\alpha+\beta}
C_{-\beta,\alpha+\beta}^{\alpha}
\\
\sum_{\alpha,\beta} (-1)^{|\alpha||\beta|}
A_{\alpha}\phi_{\beta}
C_{\alpha,\beta}^{\alpha+\beta}e_{\alpha+\beta} &\otimes&
e_{-\alpha} \otimes e_{-\beta}   \\
&:&
(-1)^{|\beta||\alpha+\beta|} A_{-\beta} \phi_{\alpha+\beta}
C_{-\beta,\alpha+\beta}^{\alpha}
\\
\sum_{\alpha,\beta}\phi_{\alpha}A_{\beta} 
C_{-\alpha,\beta}^{-\alpha+\beta}
e_{\alpha} &\otimes& e_{-\alpha+\beta} \otimes e_{-\beta}  \\
&:&
\phi_{\alpha}A_{\alpha+\beta}
C_{-\alpha,\alpha+\beta}^{\beta}
\\
\sum_{\alpha,\beta} A_{\alpha}\phi_{\beta} 
C_{-\alpha,\beta}^{-\alpha+\beta} e_{\alpha}
&\otimes& e_{-\alpha+\beta} \otimes e_{-\beta} \\
&:&
A_{\alpha}\phi_{\alpha+\beta}
C_{-\alpha,\alpha+\beta}^{\beta}
\\
\sum_{\alpha,\beta} (-1)^{|\alpha||\beta|}
\phi_{\alpha}A_{\beta} C_{-\alpha,-\beta}^{-\alpha-\beta}
e_{\alpha} &\otimes& e_{\beta} \otimes e_{-\alpha-\beta}  \\
&:&
(-1)^{|\alpha||\beta|}
\phi_{\alpha}A_{\beta} C_{-\alpha,-\beta}^{-\alpha-\beta}
\\ 
\sum_{\alpha,\beta} (-1)^{|\alpha||\beta|}
A_{\alpha}\phi_{\beta} C_{-\alpha,-\beta}^{-\alpha-\beta}
e_{\alpha} &\otimes& e_{\beta} \otimes e_{-\alpha-\beta} \\
&:&
(-1)^{|\alpha||\beta|}
A_{\alpha}\phi_{\beta} C_{-\alpha,-\beta}^{-\alpha-\beta}
\end{eqnarray*}

Gathering like terms together we see that the coefficient we
want is: 
\begin{eqnarray*} 
&& (-1)^{|\beta||\alpha+\beta|} \phi_{-\beta}A_{\alpha+\beta}
C_{-\beta,\alpha+\beta}^{\alpha}
+
(-1)^{|\beta||\alpha+\beta|} A_{-\beta} \phi_{\alpha+\beta}
C_{-\beta,\alpha+\beta}^{\alpha}\\
&+& \phi_{\alpha}A_{\alpha+\beta}
C_{-\alpha,\alpha+\beta}^{\beta}
+ 
A_{\alpha}\phi_{\alpha+\beta}
C_{-\alpha,\alpha+\beta}^{\beta} \\
&+& (-1)^{|\alpha||\beta|}
\phi_{\alpha}A_{\beta} C_{-\alpha,-\beta}^{-\alpha-\beta}
+
(-1)^{|\alpha||\beta|}
A_{\alpha}\phi_{\beta} C_{-\alpha,-\beta}^{-\alpha-\beta} 
\end{eqnarray*}
\noindent
or equivalently:
\begin{eqnarray*}
&=&\phi_{\alpha}A_{\alpha+\beta}
C_{-\alpha,\alpha+\beta}^{\beta}
+ (-1)^{|\alpha||\beta|}
\phi_{\alpha}A_{\beta} C_{-\alpha,-\beta}^{-\alpha-\beta}
\\
&&+ (-1)^{|\beta||\alpha+\beta|} \phi_{-\beta}A_{\alpha+\beta}
C_{-\beta,\alpha+\beta}^{\alpha}
+
(-1)^{|\alpha||\beta|}
A_{\alpha}\phi_{\beta} C_{-\alpha,-\beta}^{-\alpha-\beta} 
\\
&&+ A_{\alpha}\phi_{\alpha+\beta}
C_{-\alpha,\alpha+\beta}^{\beta} 
+ (-1)^{|\beta||\alpha+\beta|} A_{-\beta} \phi_{\alpha+\beta}
C_{-\beta,\alpha+\beta}^{\alpha}
\\
&=& \phi_{\alpha}\left(
A_{\alpha+\beta}
C_{-\alpha,\alpha+\beta}^{\beta} + (-1)^{|\alpha||\beta|}
A_{\beta} C_{-\alpha,-\beta}^{-\alpha-\beta} \right) 
\\
&&+ \phi_{\beta}\left(
-(-1)^{|\beta|}(-1)^{|\beta||\alpha+\beta|}A_{\alpha+\beta}
C_{-\beta,\alpha+\beta}^{\alpha}
+ (-1)^{|\alpha||\beta|}
A_{\alpha}C_{-\alpha,-\beta}^{-\alpha-\beta} \right) 
\\
&&+ \phi_{\alpha+\beta} \left(
A_{\alpha} C_{-\alpha,\alpha+\beta}^{\beta} 
+ (-1)^{|\beta||\alpha+\beta|} A_{-\beta} 
C_{-\beta,\alpha+\beta}^{\alpha} \right)
\end{eqnarray*}

We have computed earlier that:

\begin{eqnarray*}
C_{-\beta,\alpha+\beta}^{\alpha}
&=&
\frac{(-1)^{|\alpha||\beta|}(-1)^{|\beta|}
A_{-\beta}}{C_{\alpha,\beta}^{\alpha+\beta}}
(h_{\alpha},h_{\beta}) 
\\
C_{-\alpha,\alpha+\beta}^{\beta} 
&=& 
\frac{-(-1)^{|\alpha|}A_{-\alpha}}
{C_{\alpha,\beta}^{\alpha+\beta}} (h_{\alpha},h_{\beta}) 
\\
C_{-\alpha,-\beta}^{-\alpha-\beta} 
&=& 
\frac{(-1)^{|\alpha||\beta|}(-1)^{|\alpha|+|\beta|}A_{\alpha+\beta}
A_{-\alpha}A_{-\beta}}
{C_{\alpha,\beta}^{\alpha+\beta}} (h_{\alpha},h_{\beta}) 
\end{eqnarray*}
\noindent

Using these formulas, we can now see that each
pair adds up to zero, and this completes the proof of the lemma.
$\blacksquare$

This lemma shows that in order to study zero-weight
super dynamical r-matri-ces with nonzero coupling constant, it
suffices to solve the modified dynamical Yang-Baxter equation.
This modified version of the dynamical Yang-Baxter equation: 
\begin{equation*}
\tag{\ref{MDYBequation}} Alt_s(ds) + [[s,s]] +
\frac{\epsilon^2}{4}[[\Omega,\Omega]] = 0, 
\end{equation*}
\noindent 
is an equation in $\g \otimes \g \otimes \g$. 
The non-Casimir components are easily seen to be skew-symmetric
with respect to signed permutations of factors, following
arguments from the proof of Theorem \ref{0couple0weighttheorem}.
The Casimir part is also skew-symmetric
with respect to signed permutations of factors,
as can be seen by a straight-forward computation.  
Because of this symmetry, it is sufficient to look only
at the ${\h \otimes \h \otimes \h,}$ ${\h \otimes \g_{\alpha}
\otimes \g_{-\alpha},}$ and ${\g_{\alpha} \otimes \g_{\beta}
\otimes \g_{-\alpha-\beta}}$ components in order to solve
Equation \ref{MDYBequation}.

The ${\h \otimes \h \otimes \h}$ part comes only from the
dynamical part (in other words, from ${Alt_s(ds)}$) and yields:
$$ 
\frac{\partial D_{jk}}{\partial x_i} + 
\frac{\partial D_{ki}}{\partial x_j} +
\frac{\partial D_{ij}}{\partial x_k} = 0, $$
\noindent
which implies that ${D = \sum_{i<j} D_{ij} dx_i \wedge dx_j}$ is
a closed differential $2-$form.

The ${\h \otimes \g_{\alpha} \otimes \g_{-\alpha}}$ part
consists of two components. The part that comes from the
non-Casimir part is the same as the whole of the ${\h \otimes 
\g_{\alpha} \otimes \g_{-\alpha}}$ part of the dynamical
Yang-Baxter equation for a zero-weight super dynamical
r-matrix with zero coupling constant, (see Section
\ref{Proofof00Thm}):

\begin{eqnarray*}  
&&Alt_s(ds) \left|_{\h \otimes \g_{\alpha}
\otimes \g_{-\alpha}} \right. +  \sum_{\alpha,\beta}
(-1)^{|\alpha||\beta|} \delta_{\alpha,-\beta} \phi_{\alpha}
\phi_{\beta} [e_{\alpha},e_{\beta}] \otimes e_{-\alpha} \otimes
e_{-\beta} \\ 
&=& \sum_{i,\alpha} 
\frac{\partial \phi_{\alpha}}{\partial x_i} x_i \otimes
e_{\alpha} \otimes e_{-\alpha} +
\sum_{\alpha} (-1)^{|\alpha|} 
\phi_{\alpha} \phi_{-\alpha} (-A_{\alpha}) h_{\alpha} \otimes
e_{\alpha} \otimes e_{-\alpha} \\
&=& \sum_{i,\alpha} 
\frac{\partial \phi_{\alpha}}{\partial x_i} x_i \otimes
e_{\alpha} \otimes e_{-\alpha} +
\sum_{\alpha} A_{\alpha}
\phi_{\alpha} \phi_{\alpha}  h_{\alpha} \otimes
e_{\alpha} \otimes e_{-\alpha} 
\end{eqnarray*}

To this we need to add the Casimir component, i.e. the terms
that come from ${\frac{\epsilon^2}{4}[[\Omega,\Omega]]}$:
\begin{eqnarray*}
[[\Omega,\Omega]] \left|_{\h\otimes \g_{\alpha}\otimes
\g_{-\alpha}}\right. 
&=&
\sum_{\alpha,\beta}
(-1)^{|\alpha||\beta|} \delta_{\alpha,-\beta} A_{\alpha}
A_{\beta} [e_{\alpha},e_{\beta}] \otimes e_{-\alpha} \otimes
e_{-\beta}  \\
&=&
\sum_{\alpha}
(-1)^{|\alpha|} A_{\alpha} A_{-\alpha} [e_{\alpha},e_{-\alpha}]
\otimes e_{-\alpha} \otimes e_{\alpha}  \\
&=&
\sum_{\alpha}
(-1)^{|\alpha|} A_{\alpha} A_{-\alpha} A_{-\alpha} h_{\alpha}
\otimes e_{-\alpha} \otimes e_{\alpha}  \\ 
&=& 
\sum_{\alpha} A_{-\alpha}  h_{\alpha}
\otimes e_{-\alpha} \otimes e_{\alpha}  \\ 
&=& 
-\sum_{\alpha} A_{\alpha}  h_{\alpha}
\otimes e_{\alpha} \otimes e_{-\alpha} 
\end{eqnarray*} 
\noindent
and so the ${\h \otimes \g_{\alpha} \otimes \g_{-\alpha}}$ part
of Equation \ref{MDYBequation} is:
$$\sum_{i,\alpha} \frac{\partial \phi_{\alpha}}{\partial x_i} x_i \otimes
e_{\alpha} \otimes e_{-\alpha} 
+
\sum_{\alpha} A_{\alpha} \phi_{\alpha}^2 h_{\alpha} \otimes
e_{\alpha} \otimes e_{-\alpha} 
-
\frac{\epsilon^2}{4}\sum_{\alpha} A_{\alpha}  h_{\alpha}
\otimes e_{\alpha} \otimes e_{-\alpha} $$
\noindent
which we can rewrite as:
$$\sum_{\alpha \in \Delta} \left( \sum_{i} 
\frac{\partial \phi_{\alpha}}{\partial x_i} x_i 
+ A_{\alpha} (\phi_{\alpha}^2 - \frac{\epsilon^2}{4}) h_{\alpha} 
\right) \otimes e_{\alpha} \otimes e_{-\alpha}
$$
\noindent
For this term to vanish we must have, for all ${\alpha \in
\Delta}$:  
$$ \sum_{i}  \frac{\partial \phi_{\alpha}}{\partial
x_i} x_i  +  A_{\alpha} (\phi_{\alpha}^2 -
\frac{\epsilon^2}{4}) h_{\alpha} = 0.$$ 
\noindent
We can rewrite this as: 
$$ d\phi_{\alpha} + A_{\alpha} (\phi_{\alpha}^2 -
\frac{\epsilon^2}{4}) dh_{\alpha} = 0. $$
\noindent
If we define $\mu_{\alpha}$ by: 
$$ \mu_{\alpha} = \left\{ 
\begin{array}{cl}
\sqrt{-1} & \textmd{if } \alpha \textmd{ is an odd positive root}
\\  1 & \textmd{otherwise}
\end{array}  \right.$$
\noindent
then $\mu_{\alpha}^2 = A_{\alpha}$, and the
equation we need to solve is: 
$$ d\phi_{\alpha} + \mu_{\alpha}^2 (\phi_{\alpha}^2 -
\frac{\epsilon^2}{4}) dh_{\alpha} = 0. $$

We assume $\phi_{\alpha}^2 \neq \frac{\epsilon^2}{4}$ and let
$u_{\alpha}= \mu_{\alpha}\phi_{\alpha}$. Also let
$\epsilon_{\alpha} = \frac{\mu_{\alpha}\epsilon}{2}$.
Separating variables to integrate we obtain:  
$$ \frac{1}{\mu_{\alpha}} 
\int \frac{du_{\alpha}}{\epsilon_{\alpha}^2 - u_{\alpha}^2} =
\int dh_{\alpha}  \; \; \; \Rightarrow \; \; \;  
u_{\alpha} = 
\mu_{\alpha} \phi_{\alpha} = \epsilon_{\alpha} 
\textmd{ coth }(\epsilon_{\alpha}\mu_{\alpha}h_{\alpha} +C)
$$  
\noindent
and we get: 
$$ \phi_{\alpha} = \frac{\epsilon}{2} \textmd{ coth } 
\left(\frac{A_{\alpha}\epsilon}{2}(h_{\alpha} -
\nu_{\alpha})\right)$$
\noindent
for some ${\nu_{\alpha} \in \C}$. Here $h_{\alpha}$ is viewed as
a linear function on $\h^*$ via $h_{\alpha}(\lambda) = {(\alpha,
\lambda)}$.

Finally we look at the $\g_{\alpha}
\otimes \g_{\beta} \otimes \g_{-\alpha-\beta}$ part of Equation
\ref{MDYBequation}. There is no contribution from the dynamical
part. The terms from $[[s,s]]$ are (see Section
\ref{Proofof00Thm} for more details):  
\begin{eqnarray*}
\sum_{\alpha,\beta} (-1)^{|\alpha||\beta|}
\phi_{\alpha} \phi_{\beta} [e_{\alpha},e_{\beta}] \otimes
e_{-\alpha} \otimes e_{-\beta} +
\sum_{\alpha,\beta} 
\phi_{\alpha} \phi_{\beta} e_{\alpha} \otimes
[e_{-\alpha}, e_{\beta}] \otimes e_{-\beta}
\\ +
\sum_{\alpha,\beta} (-1)^{|\alpha||\beta|}
\phi_{\alpha} \phi_{\beta} e_{\alpha} \otimes e_{\beta} \otimes
[e_{-\alpha}, e_{-\beta}]
\end{eqnarray*}
\noindent
and the coefficient for the ${e_{\alpha}\otimes
e_{\beta}\otimes e_{-\alpha-\beta}}$ term coming from this
component is:
$$(-1)^{|\beta||\alpha+\beta|} \phi_{-\beta} \phi_{\alpha+\beta}
C_{-\beta,\alpha+\beta}^{\alpha} 
+ \phi_{\alpha} \phi_{\alpha+\beta}
C_{-\alpha,\alpha+\beta}^{\beta} 
+ (-1)^{|\alpha||\beta|} \phi_{\alpha} \phi_{\beta}
C_{-\alpha,-\beta}^{-\alpha-\beta} $$
\noindent
We can simplify this further:
\[
(-1)^{|\beta||\alpha+\beta|} \phi_{-\beta} \phi_{\alpha+\beta}
C_{-\beta,\alpha+\beta}^{\alpha} 
+ \phi_{\alpha} \phi_{\alpha+\beta}
C_{-\alpha,\alpha+\beta}^{\beta} 
+ (-1)^{|\alpha||\beta|} \phi_{\alpha} \phi_{\beta}
C_{-\alpha,-\beta}^{-\alpha-\beta} 
\]
\begin{eqnarray*}
&=& 
(-1)^{|\beta||\alpha+\beta|} \phi_{-\beta} \phi_{\alpha+\beta}
\frac{(-1)^{|\alpha||\beta|}(-1)^{|\beta|}
A_{-\beta}}{C_{\alpha,\beta}^{\alpha+\beta}}
(h_{\alpha},h_{\beta}) \\
&&+ \;\;\;
\phi_{\alpha} \phi_{\alpha+\beta}
\frac{-(-1)^{|\alpha|}A_{-\alpha}}
{C_{\alpha,\beta}^{\alpha+\beta}} (h_{\alpha},h_{\beta}) \\
&&+ \;\;\; 
(-1)^{|\alpha||\beta|} \phi_{\alpha} \phi_{\beta}
\frac{(-1)^{|\alpha||\beta|}(-1)^{|\alpha|+|\beta|}A_{\alpha+\beta}
A_{-\alpha}A_{-\beta}}
{C_{\alpha,\beta}^{\alpha+\beta}} (h_{\alpha},h_{\beta}) 
\\
&=& 
-\phi_{\beta} \phi_{\alpha+\beta}
\frac{A_{\beta}}{C_{\alpha,\beta}^{\alpha+\beta}}
(h_{\alpha},h_{\beta}) 
- \phi_{\alpha} \phi_{\alpha+\beta}
\frac{A_{\alpha}} {C_{\alpha,\beta}^{\alpha+\beta}}
(h_{\alpha},h_{\beta}) \\
&& + \;\;\;  \phi_{\alpha} \phi_{\beta}
\frac{A_{\alpha+\beta} A_{\alpha}A_{\beta}}
{C_{\alpha,\beta}^{\alpha+\beta}} (h_{\alpha},h_{\beta}) 
\end{eqnarray*}

The Casimir component is: 
\begin{align*}
\frac{\epsilon^2}{4} (\sum_{\alpha,\beta} (-1)^{|\alpha||\beta|}
A_{\alpha} A_{\beta} [e_{\alpha},e_{\beta}] \otimes
e_{-\alpha} \otimes e_{-\beta} +
\sum_{\alpha,\beta} 
A_{\alpha} A_{\beta} e_{\alpha} \otimes
[e_{-\alpha}, e_{\beta}] \otimes e_{-\beta}
\\ +
\sum_{\alpha,\beta} (-1)^{|\alpha||\beta|}
A_{\alpha} A_{\beta} e_{\alpha} \otimes e_{\beta} \otimes
[e_{-\alpha}, e_{-\beta}])
\end{align*}
\noindent
and the coefficient in front of the ${e_{\alpha}\otimes
e_{\beta}\otimes e_{-\alpha-\beta}}$ term coming from this
component is:
\begin{align*}
\frac{\epsilon^2}{4}
\left ( (-1)^{|\beta||\alpha+\beta|} A_{-\beta} A_{\alpha+\beta}
C_{-\beta,\alpha+\beta}^{\alpha} \right . 
&+ A_{\alpha} A_{\alpha+\beta}
C_{-\alpha,\alpha+\beta}^{\beta} \\
& \left . + (-1)^{|\alpha||\beta|} A_{\alpha} A_{\beta}
C_{-\alpha,-\beta}^{-\alpha-\beta} \right ) 
\end{align*}
\noindent
which we can simplify as follows: 
\begin{align*}
\frac{\epsilon^2}{4} &
\left ( (-1)^{|\beta||\alpha+\beta|} A_{-\beta} A_{\alpha+\beta}
C_{-\beta,\alpha+\beta}^{\alpha} 
+ A_{\alpha} A_{\alpha+\beta}
C_{-\alpha,\alpha+\beta}^{\beta} \right . \\
& \left . + (-1)^{|\alpha||\beta|} A_{\alpha} A_{\beta}
C_{-\alpha,-\beta}^{-\alpha-\beta} \right ) \\ \\
=& \quad
\frac{\epsilon^2}{4}
\left ( (-1)^{|\alpha||\beta|} A_{\beta} A_{\alpha+\beta}
C_{-\beta,\alpha+\beta}^{\alpha} 
+ A_{\alpha} A_{\alpha+\beta}
C_{-\alpha,\alpha+\beta}^{\beta} \right . \\
& \left . + (-1)^{|\alpha||\beta|} A_{\alpha} A_{\beta}
C_{-\alpha,-\beta}^{-\alpha-\beta} \right ) \\ \\
=& \quad
\frac{\epsilon^2}{4} \left(
((-1)^{|\alpha||\beta|} A_{\beta} A_{\alpha+\beta} 
\frac{(-1)^{|\alpha||\beta|}(-1)^{|\beta|}
A_{-\beta}}{C_{\alpha,\beta}^{\alpha+\beta}}
(h_{\alpha},h_{\beta}) \right. 
\\
&+ A_{\alpha} A_{\alpha+\beta}
\frac{-(-1)^{|\alpha|}A_{-\alpha}}
{C_{\alpha,\beta}^{\alpha+\beta}} (h_{\alpha},h_{\beta}) 
\\
&+ \left. (-1)^{|\alpha||\beta|} A_{\alpha} A_{\beta}
\frac{(-1)^{|\alpha||\beta|}(-1)^{|\alpha|+|\beta|}A_{\alpha+\beta}
A_{-\alpha}A_{-\beta}}
{C_{\alpha,\beta}^{\alpha+\beta}} (h_{\alpha},h_{\beta}) \right)
\\ 
=& \quad
\frac{\epsilon^2}{4} \left(
((-1)^{|\alpha||\beta|}  A_{\alpha+\beta}
\frac{(-1)^{|\alpha||\beta|}}{C_{\alpha,\beta}^{\alpha+\beta}}
(h_{\alpha},h_{\beta}) +
A_{\alpha+\beta} \frac{-1} {C_{\alpha,\beta}^{\alpha+\beta}}
(h_{\alpha},h_{\beta}) \right. \\ 
&+ \left. (-1)^{|\alpha||\beta|} 
\frac{(-1)^{|\alpha||\beta|}A_{\alpha+\beta}}
{C_{\alpha,\beta}^{\alpha+\beta}} (h_{\alpha},h_{\beta}) \right)
\\  
=& \quad	
\frac{\epsilon^2}{4} \left(
\frac{A_{\alpha+\beta}}{C_{\alpha,\beta}^{\alpha+\beta}}
(h_{\alpha},h_{\beta}) \right)
\end{align*}

Therefore the coefficient of the ${e_{\alpha}\otimes
e_{\beta}\otimes e_{-\alpha-\beta}}$ term is:
\begin{eqnarray*}
&-&\phi_{\beta} \phi_{\alpha+\beta}
\frac{A_{\beta}}{C_{\alpha,\beta}^{\alpha+\beta}}
(h_{\alpha},h_{\beta}) 
- \phi_{\alpha} \phi_{\alpha+\beta}
\frac{A_{\alpha}} {C_{\alpha,\beta}^{\alpha+\beta}}
(h_{\alpha},h_{\beta}) \\
&+& \phi_{\alpha} \phi_{\beta}
\frac{A_{\alpha+\beta} A_{\alpha}A_{\beta}}
{C_{\alpha,\beta}^{\alpha+\beta}} (h_{\alpha},h_{\beta}) 
+
\frac{\epsilon^2}{4} \left(
\frac{A_{\alpha+\beta}}{C_{\alpha,\beta}^{\alpha+\beta}}
(h_{\alpha},h_{\beta}) \right) 
\end{eqnarray*}
\noindent
This is equal to zero if and only if:
$$\phi_{\beta} \phi_{\alpha+\beta} A_{\beta} + 
\phi_{\alpha} \phi_{\alpha+\beta} A_{\alpha} = 
\phi_{\alpha} \phi_{\beta} A_{\alpha+\beta} A_{\alpha}A_{\beta}
+ \frac{\epsilon^2}{4} A_{\alpha+\beta},$$
\noindent
or equivalently: 
\begin{equation}
\label{phirelfor0thm}
A_{\alpha+\beta} \phi_{\alpha} \phi_{\beta} + 
\frac{\epsilon^2}{4} A_{\alpha+\beta} A_{\alpha}A_{\beta}  
= \phi_{\alpha+\beta} (A_{\alpha} \phi_{\beta} + 
A_{\beta} \phi_{\alpha}) 
\end{equation}
\noindent
where we use $A_{\alpha}^2 = 1$.

Let $X = \{\alpha \in \Delta | \phi_{\alpha}^2 \neq
\frac{\epsilon^2}{4}\}$. Then clearly $X$ is closed under
changing signs because $s$ satisfies the unitarity condition.
Equation \ref{phirelfor0thm} implies that $X$ is also closed
under addition. If $\alpha, \beta \in X$ are two positive roots,
then: $$ \phi_{\alpha} = \frac{\epsilon}{2} \textmd{ coth } 
\left(\frac{A_{\alpha}\epsilon}{2}(h_{\alpha} -
\nu_{\alpha})\right)$$
$$ \phi_{\beta} = \frac{\epsilon}{2} \textmd{ coth } 
\left(\frac{A_{\beta}\epsilon}{2}(h_{\beta} -
\nu_{\beta})\right)$$
\noindent
and Equation \ref{phirelfor0thm} implies: 
$$ \phi_{\alpha+\beta} = \frac{\epsilon}{2}
A_{\alpha+\beta} \textmd{ coth } 
\left(\frac{\epsilon}{2}(h_{\alpha} - \nu_{\alpha}) +
\frac{\epsilon}{2}(h_{\beta} - \nu_{\beta}) \right)$$
\noindent
is not constant, so we must have:
$$ \phi_{\alpha+\beta} = \frac{\epsilon}{2} \textmd{ coth } 
\left(\frac{A_{\alpha+\beta}\epsilon}{2}(h_{\alpha+\beta} -
\nu_{\alpha+\beta})\right).$$
\noindent
Since $coth (-x) = -coth (x)$, this implies also that:
$$ \nu_{\alpha+\beta} = \nu_{\alpha} + \nu_{\beta}. $$
\noindent
Using $\phi_{-\alpha} = -(-1)^{|\alpha|}\phi_{\alpha}$ we can
also see that $\nu_{-\alpha} = -\nu_{\alpha}$. Therefore we can
conclude that there is some ${\nu \in \h^*}$ such that
${\nu_{\alpha} = (\alpha,\nu)}$ for all ${\alpha \in X}$.

Thus we have proved that $r$ is a super dynamical r-matrix with
zero weight and coupling constant $\epsilon$ if and only if $s =
r - \frac{\epsilon}{2}\Omega$ is of the following form:
$$ s(\lambda) = \sum_{i,j=1}^N D_{ij}(\lambda) x_i \otimes x_j +
\sum_{\alpha \in \Delta}
\phi_{\alpha} e_{\alpha} \otimes e_{-\alpha} $$ 
\noindent
and this finishes the proof of Theorem \ref{0weighttheorem}.
$\blacksquare$

\section{Examples of Super Dynamical \textit{r}-matrices}
\label{SectionSuperDynamicalExamples}

\begin{example}
Let $\g$ be a simple Lie superalgebra with nondegenerate Killing form. 
For any choice of a triangular decomposition of $\g$, the constant r-matrix:
$$ r = \frac{\epsilon}{2}\sum_i \left( x_i\otimes x_i^* \right) 
+
\epsilon\sum_{\alpha \in \Delta^+}
\left( 
e_{-\alpha} \otimes e_{\alpha} 
\right)
$$  
\noindent
and its super twist:
$$ T_s(r) = \frac{\epsilon}{2}\sum_i \left( x_i\otimes x_i^* \right) 
+
\epsilon\sum_{\alpha \in \Delta^+}
\left( 
(-1)^{|\alpha|} e_{\alpha} \otimes e_{-\alpha} 
\right)
$$  
\noindent
are both zero-weight super dynamical r-matrices with coupling
constant $\epsilon$. 
They correspond to $X = \emptyset$ and $D \equiv 0$
in Theorem \ref{0weighttheorem}. We note that these r-matrices
are in fact solutions to the classical Yang-Baxter equation.
In fact, $r$ can clearly be constructed by the main construction theorem
of \cite{Kar1} with the trivial admissible triple
${(\Gamma_1,\Gamma_2,\tau)}$ where $\Gamma_i = \emptyset$ and
$\tau = 0$.

\begin{remark}
The Cartan components of the r-matrices above involve
${x_i^* \in \h}$, the elements of the dual basis of the basis
$\{x_i\}$ for $\h$, because it may not be always possible to
choose an orthonormal basis for $\h$ in the super case.  
\end{remark}

\begin{remark}
Clearly the only solution of the classical Yang-Baxter equation 
constructible by the main construction theorem of \cite{Kar1} that we can also
construct by Theorem \ref{0weighttheorem} is $r$ given above. To reach any
other r-matrix
constructible by this method, i.e. one that involves a nontrivial
admissible triple, we would need to drop the zero weight
condition.  
\end{remark}

\end{example}

\begin{example} 
For a non-constant dynamical example, consider the
super dynamical r-matrix constructed by Theorem
\ref{0weighttheorem} using $X = \Delta$, $D \equiv 0$ and $\nu =
0$:  
$$ r(\lambda) = \frac{\epsilon}{2} \Omega + 
\sum_{\alpha \in \Delta} \frac{\epsilon}{2} 
\emph{coth}\left(
\frac{\epsilon}{2}(\alpha,\lambda)
\right) 
(-1)^{|\alpha|} (e_{\alpha}, e_{-\alpha})
\left(
e_{\alpha} \otimes e_{-\alpha} 
\right) 
$$
\noindent
For a given triangular decomposition this can be rewritten as:
$$ r(\lambda) = \frac{\epsilon}{2} \Omega + 
\sum_{\alpha \in \Delta^+} \frac{\epsilon}{2} 
\emph{coth}\left(
\frac{\epsilon}{2}(\alpha,\lambda)
\right) 
\left(
(-1)^{|\alpha|}
e_{\alpha} \otimes e_{-\alpha} 
- e_{-\alpha} \otimes e_{\alpha}
\right). 
$$
\noindent
For this triangular decomposition, if we take the limit as
${\lambda \rightarrow \infty}$ while always staying positive,
then we see that:
$$ r(\lambda) \longrightarrow
\frac{\epsilon}{2} \Omega + 
\sum_{\alpha \in \Delta^+} \frac{\epsilon}{2} 
\left(
(-1)^{|\alpha|}
e_{\alpha} \otimes e_{-\alpha}
-e_{-\alpha} \otimes e_{\alpha}
\right) $$
\noindent
which is exactly the r-matrix $T_s(r)$ of the previous example.
If we take a similar limit as ${|\lambda| \rightarrow \infty}$
while staying always negative, we see that:
$$ r(\lambda) \longrightarrow
\frac{\epsilon}{2} \Omega + 
\sum_{\alpha \in \Delta^+} \frac{\epsilon}{2} 
\left(
e_{-\alpha} \otimes e_{\alpha}
-(-1)^{|\alpha|}
e_{\alpha} \otimes e_{-\alpha}
\right) $$
\noindent
which is exactly the r-matrix $r$ of the previous example.

Hence we can conclude that this particular super dynamical
r-matrix extrapolates different solutions of the classical
Yang-Baxter equation associated to the trivial admissible
triple, labeled by different triangular decompositions. 
\end{example}

\section{Conclusion}
\label{SectionConclusion}

The two theorems proven in this paper are the
super versions of Theorem $3.2$ and Theorem $3.10$ in \cite{EV}. Their proofs were very clearly
inspired by the non-graded case. A careful eye will discern the difficulties that are peculiar to
the super case in the proofs provided here. 

For a classification of all super dynamical $r-$matrices, one needs further study. Schiffmann's
results in \cite{Schi} provide a good template. It is this author's belief that the constructive
part of Schiffmann's results will carry over to the super case with certain modifications. 
In light of recent results, one expects to see some divergence from the non-graded theory when it
comes to the classification part; in fact the example studied in \cite{Kar2} will be one of the
challenges to a direct generalization.


\begin{thebibliography}{99}





\bibitem{BD1} Belavin, A. A., Drinfeld, V. G.; ``\textit{Solutions of 
the Classical Yang-Baxter Equation and Simple Lie Algebras}", Funct. 
Anal.
Appl. \textbf{16} (1982), pp.159--180.

\bibitem{BD2} Belavin, A. A., Drinfeld, V. G.; ``\textit{Triangle
Equation and Simple Lie Algebras}", Soviet Scientific Reviews
Sect. C \textbf{4} (1984), pp.93--165.


\bibitem{Et} Etingof, P.; ``\textit{On the Dynamical Yang-Baxter
Equation}", Proceedings of the International Congress of
Mathematicians, Vol. II (Beijing, 2002), Higher Ed.
Press, 2002, pp.555--570.


\bibitem{ESS} Etingof, P., Schedler, T., Schiffmann, O.;
``\textit{Explicit quantization of dynamical $r$-matrices for finite
dimensional semisimple Lie
algebras}"; J. Amer. Math. Soc. \textbf{13} (2000), no. 3, pp.595--609. 

\bibitem{ES} Etingof, P., Schiffmann, O.; Lectures on Quantum
Groups, International Press, 1998.

\bibitem{EV} Etingof, P., Varchenko, A.; ``\textit{Geometry and
Classification of Solutions of the Classical Dynamical
Yang-Baxter Equation}", Comm. Math. Phys. \textbf{192} (1998), no.
1, pp.77--120.


\bibitem{Ge} Geer, N.; ``\textit{Etingof-Kazhdan Quantization of Lie
Superbialgebras}"; e-arXiv
preprint, {\tt arXiv:math.QA/0409563} 



\bibitem{Kar1} Karaali, G.; ``\textit{Constructing r-matrices on
Simple Lie Superalgebras}", J.
Algebra \textbf{282} (2004), no.1, pp.83--102.

\bibitem{Kar2} Karaali, G.; ``\textit{A New Lie Bialgebra Structure on
$sl(2,1)$}", submitted.


\bibitem{Schi} Schiffmann, O.; ``\textit{On Classification of
Dynamical r-matrices}", Math. Res. Lett. \textbf{5} (1998),
pp.13--30.


\end{thebibliography}
\end{document}